\documentclass[11pt, a4paper]{amsart}

\usepackage{amsfonts,amsmath,amssymb,amscd}  \usepackage[all]{xy}

\def\Z{{\Bbb Z}}

\theoremstyle{definition}

\usepackage{color}

\newcommand\pf{\begin{proof}}
\newcommand\epf{\end{proof}}

\usepackage{mathrsfs}

\numberwithin{equation}{section}

\title[On integer values of generating functions]
{On integer values of the generating functions for sequences given by Pell's equations
}

\author[Y~Tsuno]{Yuji~Tsuno}
\address{Yuji Tsuno: 
Iiyama high school, 
Nagano 
389-2253, Japan}
\email{tsuno.yuji@gmail.com}

\begin{document} 

\begin{abstract}

D. S. Hong and P. Pongsriiam have provided a necessary and sufficient condition for the generating function for Fibonacci numbers (resp. the Lucas numbers) to be an integer value, for rational numbers. In other words, their results relate to the integer values of the generating functions of the sequences obtained from the integer solutions of  Pell's equation $5x^{2}-y^{2}=\pm4$. If we change this  Pell's equation to another type of Pell's equation, how will their results change?
This is a natural and interesting problem. In this paper, we show that a result similar to theirs is obtained for the generating functions for sequences given by Pell's equation $x^{2}-my^{2}=\pm1 \ (m\text{ is  a non-square natural number})$. 

\end{abstract}

\maketitle



\section{Previous results and main results}\label{sec:intro}
\vspace{5mm}

The Fibonacci sequence $\{F_{n}\}_{n \ge 0}$ is defined by $F_{0}=0$, $F_{1}=1$,
\[F_{n+2}=F_{n+1}+F_{n}.\]
The generating function for the Fibonacci sequence is given by
\[F(x)=\displaystyle\sum_{n=0}^{\infty}F_{n}x=\frac{x}{1-x-x^2}.\]

The Lucas sequence $\{L_{n}\}_{n \ge 0}$ is defined by 
$L_{0}=2$, $L_{1}=1$,
\[L_{n+2}=L_{n+1}+L_{n}.\]
The generating function for the Lucas sequence is given by
\[L(x)=\displaystyle\sum_{n=0}^{\infty}L_{n}x=\frac{2-x}{1-x-x^2}.\]

For the generating functions $F(x)$ and $L(x)$, D. S. Hong [2] proved that
\[\text{if }x=\frac{F_{n}}{F_{n+1}}(n \ge 0),\text{ then }F(x) \in \Z,\]
\[\text{ if }x=\frac{F_{n}}{F_{n+1}}(n \ge 0)\text{ or }x=\frac{L_{n}}{L_{n+1}}(n \ge 0),\text{ then }L(x) \in \Z.\]
Moreover, he questioned whether both the generating functions would be integers only in these cases. To answer this question, P. Pongsriiam [3] provided a necessary and sufficient condition for the generating function for the Fibonacci numbers (resp. the Lucas numbers) to be an integer value, for rational numbers.

\vspace{3mm}

\noindent{\bf Theorem 1.1.}(Hong$\cdot$Pongsriiam  2017)
Let $x$ be a rational number. For the generating function $F(x)$, we have $F(x) \in \Z$ if and only if\[x=\frac{F_{n}}{F_{n+1}} (n \ge 0)\text{ or  }x=-\frac{F_{n+1}}{F_{n}} (n \ge 1).\]

\vspace{3mm}

\noindent{\bf Theorem 1.2.}(Hong$\cdot$Pongsriiam 2017)
Let $x$ be a rational number. For the generating function $L(x)$, we have $L(x) \in \Z$ if and only if \[x=\frac{F_{n}}{F_{n+1}},-\frac{L_{n+1}}{L_{n}} , \frac{L_{n}}{L_{n+1}} (n \ge 0)\text{ or }x=-\frac{F_{n+1}}{F_{n}} (n \ge 1).\]

\vspace{3mm}

In the proof of the above theorems, the following famous equations for Fibonacci and Lucas numbers were used:

\[F_{n-1}F_{n+1}-F_{n}^{2}=(-1)^{n} \ (n \ge 1),\]
\[L_{n}F_{m}=F_{n+m}-(-1)^{m}F_{n-m} \  (n \ge m), \] \[5F_{n}F_{m}=L_{n+m}-(-1)^{m}L_{n-m} \ (n \ge m),\]
\[F_{n}L_{m}-L_{n}F_{m}=2(-1)^{m}F_{n-m} \ (n \ge m).\]

\vspace{3mm}

Furthermore, the following theorem for Fibonacci and Lucas numbers is important in the proof of the above theorems.

\vspace{3mm}

\noindent{\bf Theorem 1.3.}
Let $(x,y)$ be a pair of non-negative integers.
If $(x,y)$ satisfies Pell's equation
\[5x^{2}-y^{2}=\pm4,\] 
there exists a non-negative integer $n$ such that $x=F_{n}, y=L_{n}$. Conversely, 
\[5F_{n}^{2}-L_{n}^{2}=4(-1)^{n}\] 
 for any non-negative integer $n$.

\vspace{3mm}
Pell's equation \[5x^{2}-y^{2}=\pm4\] is important in the previous studies discussed. In other words, their results relate to the integer values of the generating functions of the sequences obtained from the integer solutions of  Pell's equation $5x^{2}-y^{2}=\pm4$. If we change this  Pell's equation to another type of Pell's equation, how will their results change? This is a natural and interesting problem. In this paper,
 we consider the Pell's equation 
\[x^{2}-my^{2}=\pm1 \ (m\text{ is  a non-square natural number}).\]
First, let $(a,b)$ be the minimum solution of Pell's equation $x^{2}-my^{2}=\pm1$. Here, $(a,b)$ is the integer solution such that $a \le x$ and $b \le y$ for any positive integer solution $(x,y)$.

Put
\[x_{n}= \frac{(a+b\sqrt{m})^{n}+(a-b\sqrt{m})^{n}}{2} \eqno{(1)}\]
\[y_{n}= \frac{(a+b\sqrt{m})^{n}-(a-b\sqrt{m})^{n}}{2\sqrt{m}} \eqno{(2)}\]
for $n \geq 0$.

Then $(x_{n}, y_{n}) (n \geq 0)$ is a solution of Pell's equation $x^{2}-my^{2}=\pm1$. Moreover, all non-negative integer solutions are given by these. (For example, [4, p214, Theorem 3.8.]) Here, if $a^{2}-mb^{2}=1$, we define the two sequences $\{L_{n}^{+}\}_{n \ge 0}$ and $\{F_{n}^{+}\}_{n \ge 0}$ by
\[x_{n}=L_{n}^{+}, \ y_{n}=F_{n}^{+}.\]
On the other hand, if $a^{2}-mb^{2}=-1$, we define the two sequences $\{L_{n}^{-}\}_{n \ge 0}$ and $\{F_{n}^{-}\}_{n \ge 0}$ by
\[x_{n}=L_{n}^{-}, \ y_{n}=F_{n}^{-}.\]
Thus, the generating function of the sequence $\{L_{n}^{+}\}_{n \ge 0}$ is given by
\[L^{+}(x)=\frac{1-ax}{1-2ax+x^{2}},\]
the generating function of the sequence
$\{L_{n}^{-}\}_{n \ge 0}$ is given by
\[L^{-}(x)=\frac{1-ax}{1-2ax-x^{2}},\]
the generating function of the sequence $\{F_{n}^{+}\}_{n \ge 0}$ is given by
\[F^{+}(x)=\frac{bx}{1-2ax+x^{2}}\] and
the generating function of the sequence $\{F_{n}^{-}\}_{n \ge 0}$ is given by
\[F^{-}(x)=\frac{bx}{1-2ax-x^{2}}.\]
These generating functions are drawn from the following equations obtained from $(1)$ and $(2)$:
\[F_{n+2}^{+}=2aF_{n+1}^{+}-F_{n}^{+} \eqno{(3)}\]
\[F_{n+2}^{-}=2aF_{n+1}^{-}+F_{n}^{-} \eqno{(4)}\]
\[L_{n+2}^{+}=2aL_{n+1}^{+}-L_{n}^{+} \eqno{(5)} \]
\[L_{n+2}^{-}=2aL_{n+1}^{-}+L_{n}^{-} \eqno{(6)}. \]
Furthermore, the convergence radii of these generating functions are all 
\[\frac{1}{a+b\sqrt{m}}.\] 

From here, we describe the main results of this study.

\vspace{3mm}

\noindent{\bf Theorem 1.4.}
Let $x$ be a rational number. We have 
$L^{+}(x) \in \Z$ if and only if
\[x=\frac{F_{n}^{+}}{F_{n+1}^{+}}, \frac{L_{n}^{+}}{L_{n+1}^{+}}, \frac{L_{n+1}^{+}}{L_{n}^{+}} \ (n \ge 0)\] or \[x=\frac{F_{n+1}^{+}}{F_{n}^{+}} \ (n \ge 1).\]

\vspace{3mm}

\noindent{\bf Theorem 1.5.}
Let $x$ be a rational number. We have $L^{-}(x) \in \Z$ if and only if
\[x=\frac{F_{n}^{-}}{F_{n+1}^{-}}, \frac{L_{n}^{-}}{L_{n+1}^{-}}, -\frac{L_{n+1}^{-}}{L_{n}^{-}} \ (n \ge 0)\] or \[x=-\frac{F_{n+1}^{-}}{F_{n}^{-}} \ (n \ge 1).\]

\vspace{3mm}

\noindent{\bf Theorem 1.6.}
Let $x$ be a rational number. We have $F^{+}(x) \in \Z$ if and only if
\[x=\frac{F_{n}^{+}}{F_{n+1}^{+}} \ (n \ge 0)\]or \[x=\frac{F_{n+1}^{+}}{F_{n}^{+}} \ (n \ge 1).\]

\vspace{3mm}

\noindent{\bf Theorem 1.7.}
Let $x$ be a rational number. We have $F^{-}(x) \in \Z$ if and only if 
\[x=\frac{F_{n}^{-}}{F_{n+1}^{-}} \ (n \ge 0)\] or \[x=-\frac{F_{n+1}^{-}}{F_{n}^{-}} \ (n \ge 1).\]

\vspace{3mm}

It is interesting that the main results have the same form as theorems given by D. S. Hong and P. Pongsriiam. These results pose the question: "Is the same result for other types of Pell's equations?"

\vspace{3mm}

We have the following corollaries from the main results:

\vspace{3mm}

\noindent{\bf Corollary 1.8.}
Let $x$ be a rational number. We assume that $x$ is in the convergence area of the generating function $L^{+}(x)$. Then, we have 
$L^{+}(x) \in \Z$ if and only if
\[x=\frac{F_{n}^{+}}{F_{n+1}^{+}}\ (n \ge 0).\]

\vspace{3mm}

\noindent{\bf Corollary 1.9.}
Let $x$ be a rational number. We assume that $x$ is in the convergence area of the generating function $L^{-}(x)$. Then, we have $L^{-}(x) \in \Z$ if and only if
\[x=\frac{F_{2n}^{-}}{F_{2n+1}^{-}}, \frac{L_{2n+1}^{-}}{L_{2n+2}^{-}} \ (n \ge 0).\]

\vspace{3mm}

\noindent{\bf Corollary 1.10.}
Let $x$ be a rational number. We assume that $x$ is in the convergence area of the generating function $F^{+}(x)$. Then, we have $F^{+}(x) \in \Z$ if and only if
\[x=\frac{F_{n}^{+}}{F_{n+1}^{+}} \ (n \ge 0).\]

\vspace{3mm}

\noindent{\bf Corollary 1.11.}
Let $x$ be a rational number. We assume that $x$ is in the convergence area of the generating function $F^{-}(x)$. Then, we have $F^{-}(x) \in \Z$ if and only if
\[x=\frac{F_{2n}^{-}}{F_{2n+1}^{-}} \ (n \ge 0).\]

\vspace{3mm}

These can be seen from the fact that $\{L_{n}^{+}\}_{n \ge 0}$, $\{L_{n}^{-}\}_{n \ge 0}$, $\{F_{n}^{+}\}_{n \ge 0}$, and $\{F_{n}^{-}\}_{n \ge 0}$ are monotonically increasing sequences and $(1), (2)$.

\vspace{3mm}

\noindent{\bf Remark 1.12.}
Let $a,b$ be positive integers. We assume that $b$ divides $a$. We consider the sequence $\{R_{n}\}_{n \ge 0}$ is defined by $R_{0}=0$, $R_{1}=1$, 
\[R_{n+2}=aR_{n+1}+bR_{n}.\]
A. Bulawa and W. K. Lee [1] provided a necessary and sufficient condition for the generating function for the sequence $\{R_{n}\}_{n \ge 0}$ to be an integer value, for rational numbers in the convergence area of this generating function. In this paper, conditionally, we mainly consider situations in which the initial conditions are different from their case.

\vspace{3mm}

\section{Proofs of the main results}
Before we prove the main results, we prepare the following equations:

\[F_{n-1}^{\pm}F_{n+1}^{\pm}-(F_{n}^{\pm})^{2}=\frac{-(\pm1)^{n-1}(a^{2}+mb^{2})+(\pm1)^{n}}{2m} \ \ (n \ge 1) \eqno{(7)}\]

\[2L_{n-1}^{\pm}L_{n+1}^{\pm}=L_{2n}^{\pm}+(\pm1)^{n-1}L_{2}^{\pm}
 \ \ (n \ge 1) \eqno{(8)}\]

\[F_{n}^{\pm}L_{m}^{\pm}-L_{n}^{\pm}F_{m}^{\pm}=(\pm1)^{m}F_{n-m}^{\pm} \ \ (n \ge m) \eqno{(9)}\]

\[L_{n}^{\pm}F_{m}^{\pm}=\frac{F_{n+m}^{\pm}-(\pm1)^{m}F_{n-m}^{\pm}}{2} \ \ (n \ge m) \eqno{(10)}\]

\[F_{n}^{\pm}F_{m}^{\pm}=\frac{L_{n+m}^{\pm}-(\pm1)^{m}L_{n-m}^{\pm}}{2m} \ \ (n \ge m) \eqno{(11)}\]

\[L_{n+1}^{\pm}=aL_{n}^{\pm}+mbF_{n}^{\pm} \ \ (n \ge  0) \eqno{(12)}\]

\[F_{n+1}^{\pm}=aF_{n}^{\pm}+bL_{n}^{\pm} \ \ (n \ge 0)  \eqno{(13)}\]

\[2(L_{n}^{\pm})^{2}=L_{2n}^{\pm}+(\pm1)^{n} \ \ (n \ge 0) \eqno{(14)}\]

\[2L_{n}^{\pm}L_{n+1}^{\pm}=L_{2n+1}^{\pm}+(\pm1)^{n}a \ \ (n \ge 0) \eqno{(15)}\]

These equations are obtained from (1) and (2).

\subsection{Proof of Theorem 1.4.}First, we show that 
\[L^{+}(\frac{F_{n}^{+}}{F_{n+1}^{+}}), L^{+}(\frac{L_{n}^{+}}{L_{n+1}^{+}}), L^{+}(\frac{L_{n+1}^{+}}{L_{n}^{+}}) \ (n \ge 0),\]
and
\[ L^{+}(\frac{F_{n+1}^{+}}{F_{n}^{+}}) \ (n \ge 1)\]
are integers.

If $n=0$, it is clear. If $n \ge 1$, using $(3)$, $(7)$ and $a^{2}-mb^{2}=1$, we obtain the following:

\begin{align*}
L^{+}(\frac{F_{n}^{+}}{F_{n+1}^{+}})&=\frac{F_{n+1}^{+}(F_{n+1}^{+}-aF_{n}^{+})}{F_{n+1}^{+}(F_{n+1}^{+}-2aF_{n}^{+})+{F_{n}^{+}}^{2}}\\
&\overset{(3)}{=}\frac{F_{n+1}^{+}(F_{n+1}^{+}-aF_{n}^{+})}{-F_{n+1}^{+}F_{n-1}^{+}+F_{n}^{2}}\\
&\overset{(7)}{=}\frac{F_{n+1}^{+}(F_{n+1}^{+}-aF_{n}^{+})}{b^{2}}
\end{align*}
Moreover, since $F_{0}^{+}=0, F_{1}^{+}=b, F_{n+2}^{+}=2aF_{n+1}^{+}-F_{n}^{+}$, $F_{n}^{+} \ (n \ge 0)$ is divided by $b$. Therefore, $L^{+}(\frac{F_{n}^{+}}{F_{n+1}^{+}}) \in \Z$ .
In the same way, we have the following:
\begin{align*}
L^{+}(\frac{F_{n+1}^{+}}{F_{n}^{+}})&=\frac{F_{n}^{+}(F_{n}^{+}-aF_{n+1}^{+})}{F_{n+1}^{+}(F_{n+1}^{+}-2aF_{n}^{+})+{F_{n}^{+}}^2}\\
&\overset{(3)}{=}\frac{F_{n}^{+}(F_{n}^{+}-aF_{n+1}^{+})}{{F_{n}^{+}}^{2}-F_{n+1}^{+}F_{n-1}^{+}}\\
&\overset{(7)}{=}\frac{F_{n}^{+}(F_{n}^{+}-aF_{n+1}^{+})}{b^{2}}
\end{align*}by using $(3)$ and $(7)$.
Therefore, $L^{+}(\frac{F_{n+1}^{+}}{F_{n}^{+}}) \in \Z$.

By using $(5)$, $(8)$, $(12)$, and $(14)$, 
\begin{align*}
L^{+}(\frac{L_{n}^{+}}{L_{n+1}^{+}})&=\frac{L_{n+1}^{+}(L_{n+1}^{+}-aL_{n}^{+})}{L_{n+1}^{+}(L_{n+1}^{+}-2aL_{n}^{+})+{L_{n}^{+}}^2}\\
&\overset{(5)}{=}\frac{L_{n+1}^{+}(L_{n+1}^{+}-aL_{n}^{+})}{-L_{n+1}^{+}L_{n-1}^{+}+{L_{n}^{+}}^2}\\
&\overset{(8)(14)}{=}\frac{L_{n+1}^{+}(L_{n+1}^{+}-aL_{n}^{+})}{-mb^2}\\
&\overset{(12)}{=}\frac{-L_{n+1}^{+}F_{n}^{+}}{b}\\
\end{align*}
Therefore, $L^{+}(\frac{L_{n}^{+}}{L_{n+1}^{+}})\in \Z$.

In the same way, we have the following:
\begin{align*}
L^{+}(\frac{L_{n+1}^{+}}{L_{n}^{+}})&=\frac{L_{n}^{+}(L_{n}^{+}-aL_{n+1}^{+})}{L_{n+1}^{+}(L_{n+1}^{+}-2aL_{n}^{+})+{L_{n}^{+}}^2}\\
&\overset{(5)}{=}\frac{L_{n}^{+}(L_{n}^{+}-aL_{n+1}^{+})}{-L_{n+1}^{+}L_{n-1}^{+}+{L_{n}^{+}}^2}\\
&\overset{(8)(14)}{=}\frac{L_{n}^{+}(L_{n}^{+}-aL_{n+1}^{+})}{-mb^2}\\
&\overset{(12)}{=}\frac{L_{n}^{+}(mb^{2}L_{n}^{+}+abmF_{n}^{+})}{mb^{2}}
\end{align*}by $(5)$,$(8)$,$(12)$, and $(14)$.

Hence, 
$L^{+}(\frac{L_{n+1}^{+}}{L_{n}^{+}}) \in \Z$ since $b$ divides $F_{n}^{+} \ (n \ge 0)$

Next, if $L^{+}(x)=k$ ($k$ is an integer) for some rational number $x$,
we show that
\[x=\frac{F_{n}^{+}}{F_{n+1}^{+}}, \frac{L_{n}^{+}}{L_{n+1^{+}}}, \frac{L_{n+1}^{+}}{L_{n}^{+}} \ (n \ge 0)\]
or 
\[x=\frac{F_{n+1}^{+}}{F_{n}^{+}} \ (n \ge 1).\]

If $k=0$, then
\[\frac{1-ax}{1-2ax+x^{2}}=0.\]
Hence, 
\[x=\frac{1}{a}=\frac{L_{0}^{+}}{L_{1}^{+}}.\]
If $k \neq 0$, then
\[\frac{1-ax}{1-2ax+x^{2}}=k.\]
Therefore,
\[kx^{2}+a(1-2k)x+k-1=0.\]
Hence,
\[x=\frac{-a(1-2k)\pm \sqrt{a^{2}(1-2k)^{2}-4k(k-1)}}{2k}.\]
Here, there exists a non-negative integer, $M$, such that
\[a^{2}(1-2k)^{2}-4k(k-1)=M^{2}\]
since $x$ is a rational number.
Moreover, 
\[M^{2}-mb^{2}(2k-1)^{2}=1\]since $a^{2}-mb^{2}=1$. Using (10), we have $F_{2N}^{+}=2L_{N}^{+}F_{N}^{+}$ for any non-negative integer $N$.
Therefore, there exists a non-negative integer $n$ such that $M=L_{2n+1}^{+}$.  
Moreover, we obtain $b(2k-1)=F_{2n+1}^{+} \ (n \ge 0)$ or $b(2k-1)=-F_{2n+1}^{+} \ (n \ge 1)$.
Hence, \[k=\frac{F_{2n+1}^{+}+b}{2b} (n \ge 0)\]
or
\[k=\frac{-F_{2n+1}^{+}+b}{2b} (n \ge 1).\]
From the above,
we have 
\[x=\frac{aF_{2n+1}^{+}+ bL_{2n+1}^{+}}{F_{2n+1}^{+}+b} \ (n \ge 0) \cdots (A),\]
\[x=\frac{aF_{2n+1}^{+}- bL_{2n+1}^{+}}{F_{2n+1}^{+}+b} \ (n \ge 0) \cdots (B),\]
\[x= \frac{-aF_{2n+1}^{+}+ bL_{2n+1}^{+}}{-F_{2n+1}^{+}+b} \ (n \ge 1) \cdots (C),\]
or
\[x= \frac{-aF_{2n+1}^{+}- bL_{2n+1}^{+}}{-F_{2n+1}^{+}+b} \ (n \ge 1)\cdots(D).\]

By transforming $ (A) $ to $ (D) $ using the equations from (7) to (15) from here,
we obtain 
\[x=\frac{F_{n}^{+}}{F_{n+1}^{+}}, \frac{L_{n}^{+}}{L_{n+1}^{+}}, \frac{L_{n+1}^{+}}{L_{n}^{+}} \ (n \ge 0)\]
or
\[x=\frac{F_{n+1}^{+}}{F_{n}^{+}} \ (n \ge 1).\]

Indeed, by transforming $(A)$, 
\begin{align*}
x&=\frac{aF_{2n+1}^{+}+ bL_{2n+1}^{+}}{F_{2n+1}^{+}+b}\\
&\overset{(10)}{=}\frac{2aL_{n+1}^{+}F_{n}^{+}+ab+ bL_{2n+1}^{+} }{2L_{n+1}^{+}F_{n}^{+}+2b}\\
&\overset{(9)(15)}{=}\frac{2aL_{n+1}^{+}F_{n}^{+}+2bL_{n}^{+}L_{n+1}^{+}}{2F_{n+1}^{+}L_{n}^{+}}\\
&\overset{(13)}{=}\frac{L_{n+1}^{+}}{L_{n}^{+}}
\end{align*}
By transforming $(B)$, 
\begin{align*}
x&=\frac{aF_{2n+1}^{+}- bL_{2n+1}^{+}}{F_{2n+1}^{+}+b}\\
&\overset{(10)}{=}\frac{2aL_{n+1}^{+}F_{n}^{+}+ab-bL_{2n+1}^{+} }{2L_{n+1}^{+}F_{n}^{+}+2b}\\
&\overset{(9)(15)}{=}\frac{2aL_{n}^{+}F_{n+1}^{+}-2bL_{n}^{+}L_{n+1}^{+}}{2F_{n+1}^{+}L_{n}^{+}}\\
&\overset{(12)(13)}{=}\frac{F_{n}^{+}}{F_{n+1}^{+}}
\end{align*}
By transforming $(C)$, 
\begin{align*}
x&=\frac{-aF_{2n+1}^{+}+ bL_{2n+1}^{+}}{
-F_{2n+1}^{+}+b}\\
&\overset{(10)}{=}\frac{-2aL_{n+1}^{+}F_{n}^{+}-ab+bL_{2n+1}^{+} }{-2L_{n+1}^{+}F_{n}^{+}}\\
&\overset{(9)(15)}{=}\frac{-2aL_{n}^{+}F_{n+1}^{+}+2bL_{n}^{+}L_{n+1}^{+}}{-2F_{n}^{+}L_{n+1}^{+}}\\
&\overset{(12)(13)}{=}\frac{L_{n}^{+}}{L_{n+1}^{+}}
\end{align*}
By transforming $(D)$, 
\begin{align*}
x&=\frac{-aF_{2n+1}^{+}- bL_{2n+1}^{+}}{
-F_{2n+1}^{+}+b}\\
&\overset{(10)}{=}\frac{-2aL_{n+1}^{+}F_{n}^{+}-ab-bL_{2n+1}^{+} }{-2L_{n+1}^{+}F_{n}^{+}}\\
&\overset{(15)}{=}\frac{-2aL_{n+1}^{+}F_{n}^{+}-2bL_{n}^{+}L_{n+1}^{+}}{-2F_{n}^{+}L_{n+1}^{+}}\\
&\overset{(13)}{=}\frac{F_{n+1}^{+}}{F_{n}^{+}}
\end{align*}

\qed

\vspace{3mm}

\subsection{Proof of Theorem 1.5.}First, we show that 
\[L^{-}(\frac{F_{n}^{-}}{F_{n+1}^{-}}), L^{-}(\frac{L_{n}^{-}}{L_{n+1}^{-}}), L^{-}(-\frac{L_{n+1}^{-}}{L_{n}^{-}}), \ (n \ge 0)\]
and
\[ L^{-}(-\frac{F_{n+1}^{-}}{F_{n}^{-}}) \ (n \ge 1)\]
are integers.

If $n =0$, it is clear.

If  $n \ge 1$, using $(4)$, $(7)$ and $a^{2}-mb^{2}=-1$,
we obtain that
\[L^{-}(\frac{F_{n}^{-}}{F_{n+1}^{-}})=\frac{F_{n+1}^{-}(F_{n+1}^{-}-aF_{n}^{-})}{(-1)^{n}b^{2}}\]
\[L^{-}(-\frac{F_{n +1}^{-}}{F_{n}^{-}})=\frac{F_{n}^{-}(F_{n}^{-}+aF_{n+1}^{-})}{(-1)^{n-1}b^{2}}\]

By using $(6)$, $(8)$, $(12)$, and $(14)$, we obtain that 
\[L^{-}(\frac{L_{n}^{-}}{L_{n+1}^{-}})=\frac{2bmL_{n+1}^{-}F_{n}^{-}}{(-1)^{n-1}2mb^{2}}\]

\[L^{-}(-\frac{L_{n+1}^{-}}{L_{n}^{-}})=\frac{2L_{n}^{-}(mb^{2}L_{n}^{-}+ambF_{n}^{-})}{(-1)^{n}2mb^{2}}\]
 $F_{n}^{-} (n\ge0)$ is divided by $b$. Therefore,
$L^{-}(\frac{F_{n}^{-}}{F_{n+1}^{-}})$, $L^{-}(-\frac{F_{n +1}^{-}}{F_{n}^{-}})$, $L^{-}(\frac{L_{n}^{-}}{L_{n+1}^{-}})$ and $L^{-}(-\frac{L_{n+1}^{-}}{L_{n}^{-}})$ are integers. 

Next, if $L^{-}(x)=k$ ($k$ is an integer) for some rational number $x$, we show that
\[x=\frac{F_{n}^{-}}{F_{n+1}^{-}}, \frac{L_{n}^{-}}{L_{n+1}^{-}}, -\frac{L_{n+1}^{-}}{L_{n}^{-}} \ (n \ge 0)\]
or
\[x=-\frac{F_{n+1}^{-}}{F_{n}^{-}} \ (n \ge 1)\].

If $k=0$, then
\[\frac{1-ax}{1-2ax-x^{2}}=0.\]
Hence, 
\[x=\frac{1}{a}=\frac{L_{0}^{-}}{L_{1}^{-}}.\]
If $k \neq 0$, then
\[\frac{1-ax}{1-2ax-x^{2}}=k\]
Hence,
\[-kx^{2}+a(1-2k)x+k-1=0.\]
Therefore,
\[x=\frac{-a(1-2k)\pm \sqrt{a^{2}(1-2k)^{2}-4k(k-1)}}{-2k}.\]
Here, since $x$ is a rational number, there exists a non-negative integer $M$ such that
\[a^{2}(1-2k)^{2}+4k(k-1)=M^{2}.\] 
Moreover, 
\[M^{2}-mb^{2}(2k-1)^{2}=-1,\]since $a^{2}-mb^{2}=-1$. For any non-negative integer $N$, we have $(L_{2N}^{-})^{2}-m(F_{2N}^{-})^{2}\neq -1$. Therefore, there exists a non-negative integer $n$ such that
$M=L_{2n+1}^{-}$.Besides,
$b(2k-1)=F_{2n+1}^{-} \ (n \ge 0)$ or $b(2k-1)=-F_{2n+1}^{-} \ (n \ge 1)$.
Hence, \[k=\frac{F_{2n+1}^{-}+b}{2b} (n \ge 0)\]
or
\[k=\frac{-F_{2n+1}^{-}+b}{2b} (n \ge 1).\]
From the above,

\[x=\frac{aF_{2n+1}^{-}+ bL_{2n+1}^{-}}{-F_{2n+1}^{-}-b} \ (n \ge 0) \cdots (E),\]
\[x=\frac{aF_{2n+1}^{-}- bL_{2n+1}^{-}}{F_{2n+1}^{-}-b} \ (n \ge 0) \cdots (F),\]
\[x= \frac{-aF_{2n+1}^{-}+ bL_{2n+1}^{-}}{F_{2n+1}^{-}-b} \ (n \ge 1) \cdots (G),\]
or
\[x= \frac{-aF_{2n+1}^{+}- bL_{2n+1}^{+}}{F_{2n+1}^{-}-b} \ (n \ge 1) \cdots (H).\]

By transforming $ (E) $ to $ (H) $ using the equations from (7) to (15) from here, we obtain 

\[x=\frac{F_{n}^{-}}{F_{n+1}^{-}}, \frac{L_{n}^{-}}{L_{n+1}^{-}}, -\frac{L_{n+1}^{-}}{L_{n}^{-}} \ (n \ge 0)\]
or
\[x=-\frac{F_{n+1}^{-}}{F_{n}^{-}} \ (n \ge 1).\]

Indeed, by transforming $(E)$,
\begin{align*}
x&=\frac{aF_{2n+1}^{-}+bL_{2n+1}^{-}}{-F_{2n+1}^{-}-b}\\
&\overset{(10)}{=}\frac{2aL_{n+1}^{-}F_{n}^{-}+(-1)^{n}ab+bL_{2n+1}^{-}}{-2L_{n+1}^{-}F_{n}^{-}-(-1)^{n}b-b}
\end{align*}
If $n$ is even,
\begin{align*}
\frac{2aL_{n+1}^{-}F_{n}^{-}+(-1)^{n}ab+bL_{2n+1}^{-}}{-2L_{n+1}^{-}F_{n}^{-}-(-1)^{n}b-b}
&\overset{(15)}{=}\frac{2aL_{n+1}^{-}F_{n}^{-}+2bL_{n}^{-}L_{n+1}^{-}}{-2L_{n+1}^{-}F_{n}^{-}-2b}\\
&\overset{(9)}{=}\frac{2aL_{n+1}^{-}F_{n}^{-}+2bL_{n}^{-}L_{n+1}^{-}}{-2F_{n+1}^{-}L_{n}^{-}}\\
&=\frac{L_{n+1}^{-}(aF_{n}^{-}+bL_{n}^{-})}{-F_{n+1}^{-}L_{n}^{-}}\\
&\overset{(13)}{=}-\frac{L_{n+1}^{-}}{L_{n}^{-}}
\end{align*}

If $n$ is odd,
\begin{align*}
\frac{2aL_{n+1}^{-}F_{n}^{-}+(-1)^{n}ab+bL_{2n+1}^{-}}{-2L_{n+1}^{-}F_{n}^{-}-(-1)^{n}b-b}&\overset{(15)}{=}-\frac{aF_{n}^{-}+bL_{n}^{-}}{F_{n}^{-}}\\
&\overset{(13)}{=}-\frac{F_{n+1}^{-}}{F_{n}^{-}}
\end{align*}

By transforming $(F)$,
\[x=\frac{aF_{2n+1}^{-}-bL_{2n+1}^{-}}{-F_{2n+1}^{-}-b}
\overset{(10)(15)}{=}\frac{2aL_{n+1}^{-}F_{n}^{-}-2bL_{n}^{-}L_{n+1}^{-}+2(-1)^{n}ab}{-2L_{n+1}^{-}F_{n}^{-}-(-1)^{n}b-b}\]

If $n$ is even,
\begin{align*}
\frac{2aL_{n+1}^{-}F_{n}^{-}-2bL_{n}^{-}L_{n+1}^{-}+2(-1)^{n}ab}{-2L_{n+1}^{-}F_{n}^{-}-(-1)^{n}b-b}&=\frac{2aL_{n+1}^{-}F_{n}^{-}-2bL_{n}^{-}L_{n+1}^{-}+2ab}{-2L_{n+1}^{-}F_{n}^{-}-2b}\\
&\overset{(9)}{=}\frac{a(2L_{n+1}^{-}F_{n}^{-}+2b)-2bL_{n}^{-}L_{n+1}^{-}}{-2F_{n+1}^{-}L_{n}^{-}}\\
&\overset{(9)}=\frac{aF_{n+1}^{-}-bL_{n+1}^{-}}{F_{n+1}^{-}}\\
&\overset{(12)(13)}{=}\frac{F_{n}^{-}}{F_{n+1}^{-}}
\end{align*}
If $n$ is odd, 
\begin{align*}
\frac{2aL_{n+1}^{-}F_{n}^{-}-2bL_{n}^{-}L_{n+1}^{-}+2(-1)^{n}ab}{-2L_{n+1}^{-}F_{n}^{-}-(-1)^{n}b-b}&=\frac{2aL_{n+1}^{-}F_{n}^{-}-2bL_{n}^{-}L_{n+1}^{-}+2ab}{-2L_{n+1}^{-}F_{n}^{-}}\\
&\overset{(9)}{=}\frac{L_{n}^{-}(aF_{n+1}^{-}-bL_{n+1}^{-})}{-L_{n+1}^{-}F_{n}^{-}}\\
&\overset{(12)(13)}{=}\frac{L_{n}^{-}}{L_{n+1}^{-}}
\end{align*}

By transforming $(G)$,
\[x=\frac{-aF_{2n+1}^{-}+bL_{2n+1}^{-}}{F_{2n+1}^{-}-b}\overset{(10)(15)}{=}\frac{-2aL_{n+1}^{-}F_{n}^{-}+2bL_{n}^{-}L_{n+1}^{-}-2(-1)^{n}ab}{2L_{n+1}^{-}F_{n}^{-}+(-1)^{n}b-b}\]

If $n$ is even, 
\begin{align*}
\frac{-2aL_{n+1}^{-}F_{n}^{-}+2bL_{n}^{-}L_{n+1}^{-}-2(-1)^{n}ab}{2L_{n+1}^{-}F_{n}^{-}+(-1)^{n}b-b}&=\frac{-2aL_{n+1}^{-}F_{n}^{-}+2bL_{n}^{-}L_{n+1}^{-}-2ab}{2L_{n+1}^{-}F_{n}^{-}}\\
&\overset{(9)}{=}\frac{L_{n+1}^{-}F_{n}^{-}}{L_{n}^{-}(-aF_{n+1}^{-}+bL_{n+1}^{-})}\\
&\overset{(12)(13)}{=}\frac{L_{n}^{-}}{L_{n+1}^{-}}.
\end{align*}

If $n$ is odd,
\begin{align*}
\frac{-2aL_{n+1}^{-}F_{n}^{-}+2bL_{n}^{-}L_{n+1}^{-}-2(-1)^{n}ab}{2L_{n+1}^{-}F_{n}^{-}+(-1)^{n}b-b}&=\frac{-2aL_{n+1}^{-}F_{n}^{-}+2bL_{n}^{-}L_{n+1}^{-}+2ab}{2L_{n+1}^{-}F_{n}^{-}-2b}\\
&\overset{(9)}{=}-\frac{aF_{n+1}^{-}-bL_{n+1}^{-}}{F_{n+1}^{-}}\\
&\overset{(12)(13)}{=}\frac{F_{n}^{-}}{F_{n+1}^{-}}
\end{align*}

By transforming $(H)$,

\[x=\frac{-aF_{2n+1}^{-}-bL_{2n+1}^{-}}{F_{2n+1}^{-}-b}\overset{(10)(15)}{=}\frac{-2aL_{n+1}^{-}F_{n}^{-}-2bL_{n}^{-}L_{n+1}^{-}}{2L_{n+1}^{-}F_{n}^{-}+(-1)^{n}b-b}\]

If $n$ is even, 
\begin{align*}
\frac{-2aL_{n+1}^{-}F_{n}^{-}-2bL_{n}^{-}L_{n+1}^{-}}{2L_{n+1}^{-}F_{n}^{-}+(-1)^{n}b-b}&=\frac{-2aL_{n+1}^{-}F_{n}^{-}-2bL_{n}^{-}L_{n+1}^{-}}{2L_{n+1}^{-}F_{n}^{-}}\\
&\overset{(13)}{=}-\frac{F_{n+1}^{-}}{F_{n}^{-}}
\end{align*}

If $n$ is odd, 
\begin{align*}
\frac{-2aL_{n+1}^{-}F_{n}^{-}-2bL_{n}^{-}L_{n+1}^{-}}{2L_{n+1}^{-}F_{n}^{-}+(-1)^{n}b-b}&=\frac{-2aL_{n+1}^{-}F_{n}^{-}-2bL_{n}^{-}L_{n+1}^{-}}{2L_{n+1}^{-}F_{n}^{-}-2b}\\
&\overset{(9)}{=}-\frac{L_{n+1}^{-}(aF_{n}^{-}+bL_{n}^{-})}{F_{n+1}^{-}L_{n}^{-}}\\
&\overset{(13)}{=}-\frac{L_{n+1}^{-}}{L_{n}^{-}}
\end{align*}
\qed
\vspace{3mm}

\subsection{Proof of Theorem 1.6.} First, we show that
\[F^{+}(\frac{F_{n}^{+}}{F_{n+1}^{+}})\ (n \ge 0)\]
and
\[ F^{+}(\frac{F_{n+1}^{+}}{F_{n}^{+}}) \ (n \ge 1)\]
are integers.
If $n \ge 0$, it is clear. 
If $n \ge 1$, using $(3)$, $(7)$ and $a^{2}-mb^{2}=1$,
we obtain 
\begin{align*}
F^{+}(\frac{F_{n}^{+}}{F_{n+1}^{+}})&=\frac{bF_{n}^{+}F_{n+1}^{+}}{F_{n+1}(F_{n+1}^{+}-2aF_{n}^{+})+{F_{n}^{+}}^2}\\
&\overset{(3)}{=}\frac{bF_{n}^{+}F_{n+1}^{+}}{-F_{n+1}^{+}F_{n-1}^{+}+{F_{n}^{+}}^{2}}\\
&\overset{(7)}{=}\frac{F_{n}^{+}F_{n+1}^{+}}{b}
\end{align*}

In the same way, using $(3)$, $(7)$ and $a^{2}-mb^{2}=1$,
we obtain 
\begin{align*}
F^{+}(\frac{F_{n+1}^{+}}{F_{n}^{+}})=&=\frac{bF_{n}^{+}F_{n+1}^{+}}{F_{n+1}(F_{n+1}^{+}-2aF_{n}^{+})+{F_{n}^{+}}^2}\\
&\overset{(3)}{=}\frac{bF_{n}^{+}F_{n+1}^{+}}{-F_{n+1}^{+}F_{n-1}^{+}+{F_{n}^{+}}^{2}}\\
&\overset{(7)}{=}\frac{F_{n}^{+}F_{n+1}^{+}}{b}
\end{align*}

Since $F_{n}^{+} \ (n \ge 0)$ is divided by $b$, $F^{+}(\frac{F_{n}^{+}}{F_{n+1}^{+}})$ and $F^{+}(\frac{F_{n+1}^{+}}{F_{n}^{+}})$are integers. 

Next, if $F^{+}(x)=k$ ($k$ is an integer), for some rational number $x$,
we show that
\[x=\frac{F_{n}^{+}}{F_{n+1}^{+}}\ (n \ge 0)\]
or
\[x=\frac{F_{n+1}^{+}}{F_{n}^{+}} \ (n \ge 1).\]

If $k=0$, then
\[\frac{bx}{1-2ax+x^{2}}=0\]
Hence, 
\[x=0=\frac{F_{0}^{+}}{F_{1}^{+}}.\]
If $k \neq 0$, then
\[\frac{bx}{1-2ax+x^{2}}=k.\]
Hence,
\[kx^{2}+(-2ak-b)x+k=0.\]
Therefore,
\[x=\frac{2ak+b\pm \sqrt{(2ak+b)^{2}-4k^{2}}}{2k}.\]
Here, since $x$ is a rational number, there exists a non-negative integer $M$ such that
\[(2ak+b)^{2}-4k^{2}=M^{2}.\]
Moreover, using $a^{2}-mb^{2}=1$, we obtain 
\[(2kbm+a)^{2}-mM^{2}=1.\]
Hence, there exists a non-negative integer $n$ such that
\[L_{2n+1}^{+}=2kbm+a, F_{2n+1}^{+}=M (n \ge 0).\]
Indeed, $a\pm1$ is not divided by $mb$. If $a\pm1$ is divided by $mb$, there exists a positive integer $l$ such that
\[a=mbl\pm1.\]
But,  
\[(mbl\pm1)^{2}-mb^{2}>1\]
This contradicts.
Moreover, using (14), for any non-negative integer $N$, $L_{2N}^{+}-1$ is divided by $mb$.
Hence, $L_{2N}^{+}\pm a$ is not divided by $mb$.
Therefore, there exists a non-negative integer $n$
such that
\[L_{2n+1}^{+}=\pm(2kbm+a).\]
Furthermore, if $L_{2n+1}^{+}=-(2kbm+a)$, using (12) and (15),
$2a$ is divided by $mb$ since $L_{2n+1}^{+}-a$ is divided by $mb$. This contradicts too. 
Therefore, we obtain that 
\[L_{2n+1}^{+}=2kbm+a.\]
Hence, we obtain that
\[x=\frac{aL_{2n+1}^{+}-1+bmF_{2n+1}^{+}}{L_{2n+1}^{+}-a}(n \ge 0) \cdots (I)\]or
\[x=\frac{aL_{2n+1}^{+}-1-bmF_{2n+1}^{+}}{L_{2n+1}^{+}-a}(n \ge 0) \cdots (J).\]

By transforming $(I)$,
\[x=\frac{aL_{2n+1}^{+}-1+bmF_{2n+1}^{+}}{L_{2n+1}^{+}-a}\overset{(12)}{=}\frac{L_{2n+2}^{+}-1}{L_{2n+1}^{+}-a}
\overset{(11)}{=}\frac{F_{n+1}^{+}}{F_{n}^{+}}\]

By transforming $(J)$,
\[x=\frac{aL_{2n+1}^{+}-1-bmF_{2n+1}^{+}}{L_{2n+1}^{+}-a}\overset{(12)(13)}{=}\frac{L_{2n}^{+}-1}{L_{2n+1}^{+}-a}
\overset{(11)}{=}\frac{F_{n}^{+}}{F_{n+1}^{+}}\]

\qed
\vspace{3mm}

\subsection{Proof of Theorem 1.7.} First, we show that
\[F^{-}(\frac{F_{n}^{-}}{F_{n+1}^{-}})\ (n \ge 0)\]
and
\[ F^{-}(-\frac{F_{n+1}^{-}}{F_{n}^{-}}) \ (n \ge 1)\]
are integers.

If $n=0$, it is clear.

If $n \ge 1$, using $(4)$, $(7)$ and $a^{2}-mb^{2}=-1$,
we obtain that
\[F^{-}(\frac{F_{n}^{-}}{F_{n+1}^{-}})=\frac{F_{n}^{-}F_{n+1}^{-}}{(-1)^{n}b}\]
and
\[F^{-}(-\frac{F_{n+1}^{-}}{F_{n}^{-}})=\frac{F_{n}^{-}F_{n+1}^{-}}{(-1)^{n}b}\]

Since $F_{n}^{-} (n\ge0)$ is divided by $b$, $F^{-}(\frac{F_{n}^{-}}{F_{n+1}^{-}})$ and $F^{-}(-\frac{F_{n+1}^{-}}{F_{n}^{-}})$ are integers.

Next, if $F^{-}(x)=k$ ($k$ is an integer) for some rational number $x$, we show that
\[x=\frac{F_{n}^{-}}{F_{n+1}^{-}}\ (n \ge 0)\]
or
\[x=-\frac{F_{n+1}^{-}}{F_{n}^{-}} \ (n \ge 1).\]

If $k=0$, then
\[\frac{bx}{1-2ax-x^{2}}=0.\]
Hence,
\[x=0=\frac{F_{0}^{-}}{F_{1}^{-}}.\]
If $k \neq 0$, then
\[\frac{bx}{1-2ax-x^{2}}=k.\]
Hence,
\[-kx^{2}+(-2ak-b)x+k=0.\]
Therefore,
\[x=\frac{2ak+b\pm \sqrt{(2ak+b)^{2}+4k^{2}}}{-2k}.\]
Here, since $x$ is a rational number,
there exists a non-negative integer $M$ such that
\[(2ak+b)^{2}+4k^{2}=M^{2}.\]
Moreover, using $a^{2}-mb^{2}=-1$, we obtain 
\[(2kbm+a)^{2}-mM^{2}=-1.\]
Therefore, there exists a non-negative integer $n$ such that
\[L_{2n+1}^{-}=(-1)^{n}(2kbm+a), F_{2n+1}^{-}=M \ (n \ge 0).\]
Indeed, for any non-negative integer $N$, $(L_{2N}^{-})^{2}-m(F_{2N}^{-})^{2}\neq -1$.
Hence, there exists a non-negative integer $n$ such that
\[L_{2n+1}^{-}=\pm(2kbm+a).\]
Moreover, using (12) and (15), 
$L_{2n+1}^{-}-(-1)^{n}a$ is divided by $bm$.
Therefore, we have 
\[L_{2n+1}^{-}=(-1)^{n}(2kbm+a).\]
Hence, we have 
\[x=\frac{(-1)^{n}aL_{2n+1}^{-}+1+bmF_{2n+1}^{-}}{(-1)^{n+1}L_{2n+1}^{-}+a}\ (n \ge 1) \cdots (K)\]
and
\[x=\frac{(-1)^{n}aL_{2n+1}^{-}+1-bmF_{2n+1}^{-}}{(-1)^{n+1}L_{2n+1}^{-}+a}\ (n \ge 1) \cdots (L).\]

If $n$ is even, by transforming $(K)$,
\[x=\frac{aL_{2n+1}^{-}+bmF_{2n+1}^{-}+1}{-L_{2n+1}^{-}+a}\overset{(12)}{=}\frac{L_{2n+2}^{-}+1}{-L_{2n+1}^{-}+a}\overset{(11)}{=}-\frac{F_{n+1}^{-}}{F_{n}^{-}}.\]

If $n$ is odd, by transforming $(K)$,
\[x=\frac{-aL_{2n+1}^{-}+bmF_{2n+1}^{-}+1}{L_{2n+1}^{-}+a}\overset{(12)(13)}{=}\frac{-L_{2n}^{-}+1}{L_{2n+1}^{-}+a}\overset{(11)}{=}\frac{F_{n}^{-}}{F_{n+1}^{-}}.\]

If $n$ is even, by transforming $(L)$,
\[x=\frac{aL_{2n+1}^{-}-bmF_{2n+1}^{-}+1}{-L_{2n+1}^{-}+a}\overset{(12)(13)}{=}
\frac{L_{2n}^{-}+1}{L_{2n+1}^{-}+a}\overset{(11)}{=}\frac{F_{n}^{-}}{F_{n+1}^{-}}.\]

If $n$ is odd, by transforming $(L)$,
\[x=\frac{-aL_{2n+1}^{-}-bmF_{2n+1}^{-}+1}{L_{2n+1}^{-}+a}\overset{(12)}{=}
\frac{-L_{2n+2}^{-}+1}{L_{2n+1}^{-}+a}\overset{(11)}{=}-\frac{F_{n+2}^{-}}{F_{n+1}^{-}}.\]

\qed

\end{document}